\documentclass[12pt]{article}

\usepackage{comment}
\usepackage{amssymb}
\usepackage{epsfig}
\usepackage{rotating}
\usepackage{natbib}
\usepackage{setspace}
\usepackage{hyperref}

\def\b#1{\mbox{\boldmath $#1$}}

\newcommand{\logit}{{\rm logit}}
\newcommand{\expit}{{\rm expit}}
\newcommand{\E}{{\rm E}}
\newcommand{\se}{{\rm se}}
\newcommand{\pa}{\partial}         
\def\pr{\mbox{{\rm pr}}}

\newcommand{\tr}{^{\prime}}
\newcommand{\diag}{\mbox{diag}}
\newcommand{\al}{\alpha}
\newcommand{\be}{\beta}

\newcommand{\de}{\delta}

\newcommand{\Si}{\Sigma}
\newcommand{\la}{\lambda}
\renewcommand{\th}{\theta}

\newcommand{\convp}{\stackrel{p}\rightarrow}
\newcommand{\convd}{\stackrel{d}\rightarrow}
\newtheorem{theo}{Theorem}

\newcommand{\ind}{\:\begin{sideways} $\hspace{-.4mm}\models$
\end{sideways}\:}
\def\bl#1{\mbox{\scriptsize\boldmath {$#1$}}} 
\begin{document}

\title{Causal inference in paired two-arm experimental 
studies under non-compliance with application to prognosis 
of myocardial infarction}

\author{Francesco Bartolucci\footnote{ Department of Economics, Finance
and Statistics, University of Perugia, 06123 Perugia, Italy, {\em
email}: bart@stat.unipg.it} \and Alessio Farcomeni\footnote{Department
    of Public Health and Infectious Diseases, Sapienza - University of
  Rome, 00186 Roma, Italy, {\em email}:
  alessio.farcomeni@uniroma1.it}} 
\date{}
\maketitle\vspace*{-0,5cm}
\begin{abstract}
\begin{singlespace}
Motivated by a study about prompt coronary angiography in myocardial
infarction, we propose a method to estimate the causal effect of a
treatment in two-arm experimental studies with possible non-compliance
in both treatment and control arms. The method is based on a causal
model for repeated binary outcomes (before and after the treatment),
which includes individual covariates and latent
variables for the unobserved heterogeneity between subjects. Moreover,
given the type of non-compliance, the model assumes the existence of
three subpopulations of subjects: {\em compliers}, {\em never-takers},
and {\em always-takers}. The model is estimated by a two-step
estimator: at the first step the probability that a subject belongs to
one of the three subpopulations is estimated on the basis of the available
covariates; at the second step the causal effects are estimated
through a conditional logistic method, the implementation of which depends
on the results from the first step. Standard errors for this
estimator are computed on the basis of a sandwich formula. The
application shows that prompt coronary angiography in patients with
myocardial infarction may significantly decrease the risk of other
events within the next two years, with a log-odds of about -2. Given
that non-compliance is significant for patients being given the
treatment because of high risk conditions, classical estimators fail to
detect, or at least underestimate, this effect. \vspace*{0.5cm}

\end{singlespace}
\end{abstract}

\newpage
\section{Introduction}
It is well known that, even in experimental studies, non-compliance is a strong source of confounding in the estimation of the causal effect of a treatment, in particular when measured and/or unmeasured factors affect both the decision to comply and the reaction to the treatment. There are basically three approaches to causal inference in these circumstances. These are based on: ({\em i}) potential outcomes or counterfactuals \cite[e.g.,][]{rubin:1974,rubin:1978,holland:1986,angrist:et:al:1996,abadie:2003,rubin:2005}, ({\em ii}) marginal structural models and inverse probability estimators for these models \citep{robins:1989,robins:1994}, or ({\em iii}) directed acyclic graphs (DAGs) formalized by \cite{pearl:1995,pearl:2009}. In particular, for two-arm experimental studies with all-or-nothing compliance\footnote{all-or-nothing compliance means that the treatment may be taken or not, ruling out partial compliance; for an approach specifically tailored to partial compliance see \cite{bart:grilli:2011}}, \cite{bart:2010} developed a method that may be applied with repeated binary outcomes and is based on an modified version of the conditional logistic estimator \citep{breslow:1980,collett:1991,roth:1998,hosmer:2000}. This method is based on a DAG model with latent variables, the parameters of which have a causal interpretation. The same model may be formulated on the basis of potential outcomes. The estimator is simple to apply, but in the formulation of \cite{bart:2010} it may be applied when non-compliance is only in the treatment arm and therefore, using the terminology of \cite{angrist:et:al:1996},  there are only {\em compliers} (who always comply with the treatment) and {\it never-takers} (who never take the treatment regardless of the assigned arm).

Motivated by an original application about the effectiveness of
coronary angiography (CA) in patients with non-ST elevation acute
coronary syndrome, in this paper we extend the approach of
\cite{bart:2010} by considering cases in which non-compliance may be
also observed in the control arm. Therefore, there are three
subpopulations: {\em compliers}, {\it never-takers}, and {\it
  always-takers} (who always take the treatment regardless of the
assigned arm). In particular, we extend the causal model of
\cite{bart:2010} and, basically following the same inferential
approach, we develop a conditional likelihood estimator of the causal
effects. The latter may be simply applied. It is worth noting that these causal effects are measured on the logit scale, given that we are dealing with binary outcomes; the same scale is used in relevant approaches to causal inference \citep[e.g.,][]{tenhave:et:al:2003,vans:goet:2003,robins:rot:2004,vanderlaan:et:al:2007}. Moreover, as in \cite{bart:2010}, the adopted estimator is based on two steps. At the first step we estimate the probability that a subject is a complier, a never-taker, or an always-taker on the basis of observable covariates for this subject. At the second step, the conditional likelihood of a logistic model, based on a suitable design matrix, which is set up by using the results from the first step, is maximized by a simple Newton-Raphson algorithm. Given the two-step formulation of the estimator, we use a sandwich formula \citep{white:1982} for deriving standard errors. These may be used to test the significance of the causal parameters.

As mentioned above, we develop our methodology in connection with an
original study on CA in patients with non-ST elevation acute coronary
syndrome. In particular, we are interested in investigating whether a
prompt CA (within 48h from hospital admission) should
be recommended in light of a lower risk of recurrent cardiovascular
events after leaving the hospital. A prompt CA, together with ECG and
other exams performed on patients with coronary syndrome, may be
helpful in better calibrating an in-hospital treatment. Even if the
current guidelines of the European cardiologic society
recommend CA within 48h of hospitalization \citep{bertrand:et:al:2002},
in some hospitals patients are submitted to CA only after a few days,
or even not at all. In the cardiology literature a definite 
recommendation has not yet emerged, with some studies 
reporting equivalence of CA performed before or after 48h of hospitalization
\citep{timi:1994,boden:et:al:1998,mcc:et:al:1998},
and other studies reporting superiority of prompt CA \citep{ragmin:et:al:1999,cannon:et:al:2001,fox:et:al:2002}. 
In our data, the medium/long term effects of
coronary angiography within 48h from hospital admission have been
estimated using a control given by the usual clinical practice in the
hospital, which may or may not include the coronary
angiography; when included, the designed study planned 
to schedule it only after at least 48h from hospitalization. 
Then, subjects assigned to the treatment group are expected to
undergo CA within 48h from hospitalization, whereas patients assigned to
the control group may or may not undergo CA. When a patient in the
control group is submitted to CA, the analysis is expected to be
executed after 48h from the hospitalization. 
Patients were randomized immediately at hospitalization. 
In practice, a significant fraction of controls received CA within 48h
from hospitalization, possibly due to the need of information in order
to promptly proceed with a treatment. Furthermore, a significant fraction of
patients in the active group (treatment arm) did receive CA, but after 48h from
hospitalization, possibly due to a busy hospital
schedule which did not allow prompt CA performance. 
We consequently have a significant non-compliance in both arms,
leading to the presence of never-takers and always-takers in addition 
to compliers. Note that non-compliance 
in this example is more likely a choice of the doctor, rather than of the
patient. 

We focus on a relevant group of patients, those arriving at the
hospital with myocardial infarction. From our analyses, based on the causal inference approach here proposed, two important findings emerge. 
First of all, there is a significant {\it causal} effect of
prompt CA, with a log odds-ratio of about -2 and $p$-value equal to $0.009$. 
Hence, patients arriving
at the hospital with myocardial infarction 
should be submitted for coronary angyography within
48h, and this will help doctors in greatly decreasing the risk of recurrent events
after dismissal.
Secondly, we estimate the effects separately on
the four groups (never-takers receiving control, compliers receiving control,
compliers receiving treatment, and always-takers receiving treatment), and we observe
that the bias is
arising mostly from the always-takers. In fact, the treatment has substantially no effect on the always-takers, but we estimate a strong effect on compliers. 

The paper is organized as follows. In Section \ref{sec:data} we briefly describe
the data from the study motivating this paper. In Section \ref{sec:model} we introduce the
causal model for repeated binary response variables. The proposed two-step
estimator is described in Section \ref{sec:estimator} and its application to
the dataset deriving from the cardiology study outlined above is described
in Section \ref{sec:application}. Final conclusions are reported in Section 
\ref{sec:conclusions}. 

We implemented the estimator in an {\tt R} function that we make available to the reader upon request.
\section{Description of the Prompt Coronary Angiography data}\label{sec:data}
The multicenter trial we consider is based on the inclusion of 
patients arriving to the hospital with last episode of {\it angina pectoris} within the last 24 hours. The patients were included in the study if they were diagnosed a myocardial infarction. Patients with persistent ST elevation or who could not undergo CA were excluded from the study. 

The binary response of interest is the recurrence within 2 years
after leaving the hospital of any among: ({\em i}) another episode
of myocardial infarction, ({\em ii}) an episode of 
{\it angina pectoris} of duration 20 minutes or longer, ({\em iii}) 
other significant cardiovascular events, or ({\em iv}) death
which could be related to the current episode. The recorded data concern the
presence or absence of episodes of {\it angina pectoris}, 
myocardial infarction or other cardiovascular events  
within the last month before hospitalization, and other covariates. 
The first can be considered as a pre-treatment copy of the outcome, 
which we denote with $Y_1$. Among the covariates there are: gender, age, 
smoke, statin use, history of CHD in the family, hypertension, 
and glicemic index (GI) at hospitalization. 
We are interested in investigating the effect of a prompt CA since 
our population of patients with myocardial infarction (IMA) at
hospitalization could probably be better treated after CA, and this
could prevent further events. 

Overall, we have data on $n=1,560$ subjects, whose
characteristics are summarized as follows: there are 63\% males, 46\% smokers, 75\% have a history of CHD in the family, 31\% have hypertension, and 81\% use statines regularly.
GI has a strongly skewed distribution,
with median equal to 118 and MAD equal to 34;
moreover, the mean age is 67.5 with a standard deviation of 10.8. 

Randomization was performed with a proportion of 1:2, and in fact 66\%
of the patients are assigned to the prompt CA group. Only 52\% of the
patients actually were submitted to prompt CA. 
There was non-compliance in both groups, with more than 1/3 of the
subjects assigned to each group ending up taking the other treatment. 
More precisely, 370 subjects assigned to the prompt CA did undergo CA
later than 48h after hospitalization, and 170 patients assigned to the
control group had prompt CA. 

Given that after model selection we will conclude that GI
and use of statines are predictive of compliance (see Section \ref{sec:application}),
we study these two variables a bit more in depth here. For this aim, 
in Table \ref{compl_GI} we report the proportion of patients
belonging to the groups of not treated as assigned (assigned and received
control), always-takers (assigned to control and received
treatment), never-takers (assigned to treatment and received
control), or treated as assigned (assigned and received
treatment), given the level of GI and the use or not of Statines. 
The level of GI is discretized on the basis of the
quartiles of the empirical distribution. 
It is important to underline that the first and last groups
are made of both compliers and subjects who were by chance assigned to
the treatment they would have preferred anyway. That is, 
in the first group we have both compliers assigned to the control and
never-takers randomized to the control; in the last group we have both
compliers assigned to the treatment and always-takers who were also 
randomized to the treatment. 

\begin{table}[ht]\centering
{\small
\begin{tabular}{ccccccccccccc}
\hline\hline
&& \multicolumn4c{GI quartile} && \multicolumn4c{Use of statines}
\\\cline{3-6}\cline{8-9}
Arm & Group & 1st & 2nd & 3rd & 4th && No & Yes\\\hline
Control & Compliers + never-takers  & 0.634 & 0.702 & 0.674 & 0.638 && 0.606 & 0.676\\
 & Always-takers & 0.366 & 0.298 & 0.326 & 0.362 && 0.394 & 0.324\\\hline
Treatment & Never-takers  & 0.336 & 0.335 & 0.389 & 0.430 && 0.443 & 0.352 \\ 
                 & Compliers + always-takers & 0.664 & 0.665 & 0.611 & 0.570 && 0.557 & 0.648\\\hline
\end{tabular}}
\caption{\em Conditional proportion of the group of not treated as assigned 
(compliers + never-takers), never-takers, always-takers, 
or treated as assigned (compliers + alwyas-takers), given GI
and the use or not of statines}\label{compl_GI}
\end{table}

From the results in table Table \ref{compl_GI} it can be seen that the proportion 
of never-takers steadily increases with GI, whereas the proportion of always-takers is 
larger for the first and last quartiles. 
On the other hand, the use of statines seems to increase the compliance in both
directions, with a decrease of 7\% of always-takers and 9\% of never-takers. 
\section{The causal model}\label{sec:model}
Let $Y_1$ and $Y_2$ denote the binary outcomes of interest, let $\b V$ be a vector of observable covariates, let $Z$ be a binary variable equal to 1 when a subject is assigned to the treatment and to 0 when he/she is assigned to the control, and let $X$ be the corresponding binary variable for the treatment actually received. In the present framework $\b V$ and $Y_1$ are pre-treatment variables, whereas $Y_2$ is a post-treatment variable. Moreover, non-compliance of the subjects involved in the experimental study implies that $X$ may differ from $Z$, since we consider experimental studies in which subjects randomized to both arms can access the treatment and therefore any configuration of $(Z,X)$ may be observed.
Consequently, we assume the existence of three subpopulations of subjects enrolled in the study: {\em compliers}, {\em never-takers}, and {\em
always-takers} \citep{angrist:et:al:1996}. This rules out the presence of {\em defiers}, that is, subjects that systematically take the treatment if assigned to the control arm and vice-versa.

In the following, we introduce a latent variable model for the analysis of data deriving from the experimental study described above. This model extends that proposed by \cite{bart:2010} to deal with two-arm experimental studies of the same type in which, however, non-compliance may be only observed in the treatment arm. We then derive results about the proposed model which are useful for making inference on its parameters.
\subsection{Model assumptions}
We assume that the behavior of a subject depends on the observable covariates $\b V$, a latent variable $U$ representing the effect of unobservable covariates on both response variables, and a latent variable $C$ representing the attitude to comply with the assigned treatment. The last one, in particular, is a discrete variable with three possible values: 0 for never-takers, 1 for compliers, and 2 for always-takers. The model is based on   assumptions A1-A5 that are reported below. In formulating these assumptions we use the symbol $W_1\ind W_2|W_3$ to denote conditional independence between the random variables $W_1$ and $W_2$ given $W_3$; this notation naturally extends to random vectors. Moreover, with reference to the variables in our study, we also let $p_1(y|u,\b v)=\pr(Y_1=y|U=u,\b V=\b v)$ and $p_2(y|u,\b v,c,x)=\pr(Y_2=y|U=u,\b V=\b v,C=c,X=x)$, and by $1\{\cdot\}$ we denote the indicator function.

The model assumptions are:

\begin{itemize}
\item[A1:] $C\ind Y_1|(U,\b V)$;
\item[A2:] $Z\ind(U,Y_1,C)|\b V$;
\item[A3:] $X\ind(U,\b V,Y_1)|(C,Z)$ and, with probability 1, $X=Z$ when $C=1$ (compliers), $X=0$ when $C=0$ (never-takers), and $X=1$ when $C=2$
(always-takers);
\item[A4:] $Y_2\ind (Y_1,Z)|(U,\b V,C,X)$;
\item[A5:] for all $u,\b v,c$ and $x$, we have
\[
\logit[p_2(1|u,\b v,c,x)]-\logit[p_1(1|u,\b v)]=\b t(c,x)\tr\b\be,
\]
where 
\[
\b t(c,x) = \pmatrix{1\{c=0\}(1-x)\cr 1\{c=1\}(1-x)\cr 1\{c=1\}x\cr1\{c=2\}x}\quad\mbox{and}\quad\b\be=\pmatrix{\be_0\cr\be_1\cr\be_2\cr\be_3}.
\]
\end{itemize}

According to assumption A1 the tendency to comply depends only on $(U,\b V)$, whereas according to assumption A2 the randomization only depends on the observable covariates in $\b V$. This assumption is typically satisfied in randomized experiments of our interest and, in any case, it may be relaxed by requiring that $Z$ is conditionally independent of $U$ given $(\b V,Y_1)$; this is shown in \cite{bart:2010}.
Assumption A3 is rather obvious considering that $C$ represents the tendency of a subject to comply with the assigned treatment. Assumption A4 implies that there is no direct effect of $Y_1$ on $Y_2$, since the distribution of the latter depends only on 
$(U,\b V,C,X)$; it also implies an assumption known as {\em exclusion restriction}, according to which $Z$ affects $Y_2$ only through $X$. Finally, assumption A5 states that the distribution of $Y_2$ depends on a vector of causal parameters $\b\be$, the elements of which are interpretable as follows:
\begin{itemize}
\item $\be_0$: effect of control on never-takers;
\item $\be_1$: effect of control on compliers;
\item $\be_2$: effect of treatment on compliers;
\item $\be_3$: effect of treatment on always-takers. 
\end{itemize}

The most interesting quantity to estimate is the {\em causal effect} of the treatment over the control in the subpopulation of compliers. In the present context, this effect may be defined as
\[
\de=\logit[p_2(1|u,\b v,1,1)]-\logit[p_2(1|u,\b v,1,0)]=\be_2-\be_1
\]
and corresponds to the increase of the logit of the probability of success when $x$ goes from 0 to 1, all the other factors remaining unchanged. 

The above assumptions imply the dependence structure between the observable and unobservable variables may be represented by the same DAG reported in \cite{bart:2010}. These assumptions lead to a causal model in the sense of \cite{pearl:1995} since all the observable and unobservable factors affecting the response variables of interest are included. Moreover, using the same approach used in \cite{bart:2010}, the model may be also formulated in terms of potential outcomes, enforcing in this way its causal interpretation.
\subsection{Preliminary results}
Along the same lines as in \cite{bart:2010}, assumptions A1-A5 imply that the probability function of the conditional distribution of $(Y_1,Z,X,Y_2)$ given $(U,\b V,C)$ is equal to
\[
p(y_1,z,x,y_2|u,\b v,c)=p_1(y_1|u,\b v)q(z|v)f(x|c,z)p_2(y_2|u,\b v,c,x),
\]
where $q(z|\b v)=\pr(Z=z|\b V=\b v)$ and $f(x|c,z)=\pr(X=x|C=c,Z=z)$. After some algebra, for the conditional distribution of $(Y_1,Z,X,Y_2)$
given $(U,\b V)$ we have
\begin{equation}
p(y_1,z,x,y_2|u,\b v)=\frac{e^{(y_1+y_2)\la(u,\bl v)}}{1+e^{\la(u,\bl v)}}q(z|\b v)\sum_{c=0}^2
f(x|c,z)\frac{e^{y_2\bl t(c,x)\tr\be}}{1+e^{\la(u,\bl v)+\bl t(c,x)\tr\be}}\pi(c|u,\b v),\label{p_man}
\end{equation}
where $\la(u,\b v)=\logit[p_1(y|u,\b v)]$ and $\pi(c|u,\b v)=\pr(C=c|U=u,\b V=\b v)$. 

This probability function considerably
simplifies when $x\neq z$. In fact, for $z=1$ and $x=0$, $f(x|c,z)$
is equal to 1 when $c=0$ (never-takers) and to 0 otherwise.
Similarly, for $z=0$ and $x=1$, $f(x|c,z)$ is equal to 1 when $c=2$
(always-takers) and to 0 otherwise. We then have
\[
p(y_1,z,x,y_2|u,\b v)=\frac{e^{(y_1+y_2)\la(u,\bl v)}}{1+e^{\la(u,\bl v)}}q(z|v)
\frac{e^{y_2\bl t(c,x)\tr\bl \be}}{1+e^{\la(u,\bl v)+\bl t(c,x)\tr\bl \be}}\pi(c|u,\b v),
\]
with
\begin{equation}
c = \left\{\begin{array}{ll} 0 & \mbox{if } z=1,x=0,\cr 2 & \mbox{if
} z=0,x=1.\end{array} \right.\label{def_c}
\end{equation}
Consequently, $(Y_1,Y_2)$ is conditionally independent of $U$ given
$(\b V,Z,X,Y_+)$ and $Z\neq X$. In particular, for $Y_+=1$ we have
\[
p(y_1,y_2|u,\b v,z,x,1)=p(y_1,y_2|\b v,z,x,1)=\frac{e^{y_2\bl t(c,x)\tr\bl\be}}{1+e^{\bl t(c,x)\tr\bl\be}},
\]
with $c$ defined as in (\ref{def_c}).

When $x=z$, the conditional probability $p(y_1,z,x,y_2|u,\b v)$ has
the following expression: 
\[
p(y_1,0,0,y_2|u,\b v)=\frac{e^{(y_1+y_2)\la(u,\bl v)}}{1+e^{\la(u,\bl v)}}q(0|\b v)\sum_{c=z}^{z+1}
\frac{e^{y_2t(c,0)\tr\be}}{1+e^{\la(u,\bl v)+\bl t(c,0)\tr\bl\be}}\pi(c|u,\b v);
\]
note that sum $\sum_{c=z}^{z+1}$ is extended to $c=0,1$ for $x=z=0$ and to $c=1,2$ for $x=z=1$. The latter expression is based on a mixture between the conditional distribution of $Y_2$ for the population of compliers and that of never-takers.

Finally, consider the conditional distribution of $(Y_1,Y_2)$ given $(U,\b V,Z,X,Y_+)$, with $Y_+=Y_1+Y_2$. The probability function of
this distribution is denoted by $p(y_1,y_2|u,\b v,z,x,y_+)$ and is equal to 1 for $y_+=0,2$, whereas for $y_+=1$ it may be obtained as
\[
p(y_1,y_2|u,\b v,z,x,1)=\frac{p(y_1,z,x,y_2|u,\b v)}{p(1,z,x,0|u,\b v)+p(0,z,x,1|u,\b v)}.
\]

An interesting result deriving from (\ref{p_man}) is that, when
$x\neq z$, the latter expression does not depend on $u$ and is equal to
\[
p(y_1,y_2|\b v,z,x,1)=\frac{e^{y_2\bl t(x,x)\tr\bl\be}}{1+
e^{\bl t(x,x)\tr\bl\be}}.
\]
On the other hand, $(Y_1,Y_2)$ is no longer conditionally independent of $U$, given
$(V,Z,X, Y_+)$ and $X=Z$. However, we show below that we can approximate the corresponding conditional probability function by a function which is independent of $u$. This
is the basis for the pseudo conditional likelihood estimator of
$\b\be$ and $\de$ proposed in the next section.
\section{Pseudo conditional likelihood inference}\label{sec:estimator}
For a sample of $n$ subjects included in the two-arm experimental study, let $y_{i1}$ denote the observed value of $Y_1$ for subject $i$, $i=1,\ldots,n$,
let $y_{i2}$ denote the value of $Y_2$ for the same subject, and let $\b v_i$, $z_i$, and $x_i$ denote the corresponding values of $\b V$, $Z$, and $X$, respectively. In the following, we introduce an approach for estimating the causal parameter vector $\b\be$ which closely follows that proposed in \cite{bart:2010}. The approach relies on the maximization of a likelihood based on the probability function $p(y_1,y_2|\b v,z,x,1)$, for the cases in which $(Y_1,Y_2)$ is conditionally independent of $U$ given $(V,Z,X,Y_+)$, and on an approximated version of this function otherwise. It results a pseudo conditional likelihood estimator, in the sense of \cite{white:1982},
whose main advantage is the simplicity of use; see also \cite{bart:nigro:2012} for a related approach applied in a different field. Note that this approach requires the preliminary estimation of the probability that every subject belongs to one of the three subpopulations (compliers, never-takers, and always-takers). Overall, the approach is based on two steps that are detailed in the following.

At the {\em first step} we estimate the probabilities that a 
subject is a never-taker ($c=0$), a complier ($c=1$), or an
always-taker ($c=2$). We assume that a multinomial logit 
with the category of compliers as reference category: 
\begin{eqnarray}
\log\frac{\pi(0|\b v)}{\pi(1|\b v)}&=&\b g(\b v)\tr\b\al_0,\label{eq:ass1}\\
\log\frac{\pi(2|\b v)}{\pi(1|\b v)}&=&\b g(\b v)\tr\b\al_2.\label{eq:ass2}
\end{eqnarray}
This implies that
\begin{eqnarray}
\pi(0|\b v)&=&\frac{\exp[\b g(\b v)\tr\b\al_0]}{1+\exp[\b g(\b v)\tr\b\al_0]+\exp[\b g(\b v)\tr\b\al_2]},\label{eq:inv1}\\
\pi(1|\b v)&=&\frac{1}{1+\exp[\b g(\b v)\tr\b\al_0]+\exp[\b g(\b v)\tr\b\al_2]},\label{eq:inv2}\\
\pi(2|\b v)&=&\frac{\exp[\b g(\b v)\tr\b\al_2]}{1+\exp[\b g(\b v)\tr\b\al_0]+\exp[\b g(\b v)\tr\b\al_2]}.\label{eq:inv3}
\end{eqnarray}

Given that the assignment is randomized and does not depend on the
individual covariates, 
the parameter vectors $\b\al_0$ and $\b\al_2$ are estimated by maximizing the log-likelihood
\begin{eqnarray*}
\ell_1(\b\al_0,\b\al_2)&=&\sum_i\ell_{1i}(\b\al_0,\b\al_2),\\
\ell_{1i}(\b\al_0,\b\al_2)&=&(1-z_i)(1-x_i)\log[\pi(0|\b v_i)+\pi(1|\b v_i)]+
\sum_i(1-z_i)x_i\log\pi(2|\b v_i)\\
&+&\sum_iz_i(1-x_i)\log\pi(0|\b v_i)+
\sum_iz_ix_i\log[\pi(1|\b v_i)+\pi(2|\b v_i)].
\end{eqnarray*}
For this aim, a simple Newton-Raphson algorithm may be used, which is based on the first and second derivatives of this function. In particular, the first derivative of this function may be found as follows. First of all we write
\begin{eqnarray*}
\ell_{1i}(\b\al_0,\b\al_2)&=&(1-z_i)(1-x_i)\log\frac{\pi(0|\b v_i)+\pi(1|\b v_i)}{\pi(1|\b v_i)}+(1-z_i)x_i\log\frac{\pi(2|\b v_i)}{\pi(1|\b v_i)}\\
&+&z_i(1-x_i)\log\frac{\pi(0|\b v_i)}{\pi(1|\b v_i)}+
z_ix_i\log\frac{\pi(1|\b v_i)+\pi(2|\b v_i)}{\pi(1|\b v_i)}+n\log\pi(1|\b v_i).
\end{eqnarray*}
Then, based on the above assumptions (\ref{eq:ass1}) and (\ref{eq:ass2}), we have
\begin{eqnarray*}
\ell_{1i}(\b\al_0,\b\al_2)&=&(1-z_i)(1-x_i)\log\{1+\exp[\b g(\b v_i)\tr\b\al_0]\}+
(1-z_i)x_i\b g(\b v_i)\tr\b\al_2\\
&+&z_i(1-x_i)\b g(\b v_i)\tr\b\al_0+
z_ix_i\log\{1+\exp[\b g(\b v_i)\tr\b\al_2]\}\\
&-&\log\{1+\exp[\b g(\b v_i)\tr\b\al_0]+\exp[\b g(\b v_i)\tr\b\al_2]\},
\end{eqnarray*}
so that
\begin{eqnarray*}
\frac{\pa \ell_1(\b\al_0,\b\al_2)}{\pa\b\al_0}&=&\sum_i\frac{\pa \ell_{1i}(\b\al_0,\b\al_2)}{\pa\b\al_0},\\
\frac{\pa \ell_1(\b\al_0,\b\al_2)}{\pa\b\al_0}&=&\left[(1-z_i)(1-x_i)\pi^*(0|\b v_i)+z_i(1-x_i)-\pi(0|\b v_i)\right]\b g(\b v_i),
\end{eqnarray*}
and
\begin{eqnarray*}
\frac{\pa \ell_1(\b\al_0,\b\al_2)}{\pa\b\al_2} &=&\sum_i\frac{\pa \ell_{1i}(\b\al_0,\b\al_2)}{\pa\b\al_2},\\ 
\frac{\pa \ell_{1i}(\b\al_0,\b\al_2)}{\pa\b\al_2}&=&\left[z_ix_i\pi^*(2|\b v_i)+(1-z_i)x_i-\pi(2|\b v_i)\right]\b g(\b v_i),
\end{eqnarray*}
where
\begin{eqnarray}
\pi^*(0|\b v_i)=\frac{\pi(0|\b v_i)}{\pi(0|\b v_i)+\pi(1|\b v_i)},\quad 
\pi^*(2|\b v_i)=\frac{\pi(2|\b v_i)}{\pi(1|\b v_i)+\pi(2|\b v_i)}.\label{eq:ratio_prob}
\end{eqnarray}
Moreover, regarding the second derivative, we have
\begin{eqnarray*}
\frac{\pa^2\ell_1(\b\al_0,\b\al_2)}{\pa\b\al_0\pa\b\al_0\tr} &=& 
\sum_i\left\{(1-z_i)(1-x_i)\pi^*(0|\b v_i)[1-\pi^*(0|\b v_i)]-
\pi(0|\b v_i)[1-\pi(0|\b v_i)]\right\}\b g(\b v_i)\b g(\b v_i)\tr,\\
\frac{\pa^2\ell_1(\b\al_0,\b\al_2)}{\pa\b\al_0\pa\b\al_2\tr} &=& 
\sum_i\pi(0|\b v_i)\pi(2|\b v_i)\b g(\b v_i)\b g(\b v_i)\tr,\\
\frac{\pa^2\ell_1(\b\al_0,\b\al_2)}{\pa\b\al_0\pa\b\al_0\tr} &=& 
\sum_i\left\{z_ix_i\pi^*(2|\b v_i)[1-\pi^*(2|\b v_i)]-
\pi(2|\b v_i)[1-\pi(2|\b v_i)]\right\}\b g(\b v_i)\b g(\b v_i)\tr.
\end{eqnarray*}
The estimated parameter vectors, obtained by maximizing $\ell_1(\b\al_0,\b\al_2)$, 
are denoted by $\hat{\b\al}_0$ and $\hat{\b\al}_2$ and the corresponding 
probabilities are denoted by $\hat{\pi}(0|\b v)$,  $\hat{\pi}(1|\b v)$, and 
$\hat{\pi}(2|\b v)$, which are obtained by (\ref{eq:inv1}), (\ref{eq:inv2}), and 
(\ref{eq:inv3}), respectively. 
Finally, by inversion of minus the Hessian matrix, which is based on the second derivatives above, it is also possible to obtain the standard errors for the parameter estimates $\hat{\b\al}_0$ and $\hat{\b\al}_2$ in the usual way.

At the {\em second step}, we maximize the following weighted conditional log-likelihood:
\begin{eqnarray*}
\ell_2(\b\be|\hat{\b\al}_0,\hat{\b\al}_2) &=& \sum_id_i
\ell_{2i}(\b\be|\hat{\b\al}_0,\hat{\b\al}_2),\\
\ell_{2i}(\b\be|\hat{\b\al}_0,\hat{\b\al}_2)&=&
(1-z_i)(1-x_i)\sum_{c=0}^1\hat{\pi}_{01}^*(c|\b v_i)
\frac{\exp(y_{i2}\be_c)}{1+\exp(\be_c)}+
(1-z_i)x_i\frac{\exp(y_{i2}\be_3)}{1+\exp(\be_3)}\\
&+&z_i(1-x_i)\frac{\exp(y_{i2}\be_0)}{1+\exp(\be_0)}+
z_ix_i\sum_{c=1}^2\hat{\pi}_{12}^*(c|\b v_i)
\frac{\exp(y_{i2}\be_{c+1})}{1+\exp(\be_{c+1})},
\end{eqnarray*}
where $d_i=1\{y_{i1}+y_{i2}=1\}$, so that only discordant configurations are considered, and
\[
\eta_h = \frac{\exp(\be_h)}{1+\exp(\be_h)},\quad h=0,\ldots,3,
\]
where $\be_0$ is the effect of placebo on never-takers, $\be_1$ is the effect of placebo on compliers, $\be_2$ is the effect of treatment on compliers, and $\be_3$ is the effect of treatment of always-takers. Finally, as generalization of (\ref{eq:ratio_prob}), we have that
\[
\hat{\pi}_{01}^*(c|\b v_i)=\frac{\hat{\pi}(c|\b v_i)}{\hat{\pi}(0|\b v_i)+\hat{\pi}(1|\b v_i)},\quad c=0,1
\]
and
\[
\hat{\pi}_{12}^*(c|\b v_i)=\frac{\hat{\pi}(c|\b v_i)}{\hat{\pi}(1|\b v_i)+\hat{\pi}(2|\b v_i)},\quad c=1,2.
\]
The first is the probability of being a never-taker or a complier given that the subject is in one of these subpopulation and his/her covariates; 
a similar interpretation holds for the probabilities of the second type. 

In order to compute the first and second derivatives of 
$\ell_2(\b\be|\hat{\b\al}_0,\hat{\b\al}_2)$ with respect to $\b\be$,
it is convenient to express $i$-th component of this function as
\[
\ell_{2i}^*(\b\eta|\hat{\b\al}_0,\hat{\b\al}_2)=y_{i2}\log(\hat{\b w}_i\tr\b\eta)+
(1-y_{i2})\log(1-\hat{\b w}_i\tr\b\eta),
\]
where $\b\eta=(\eta_0,\eta_1,\eta_2,\eta_3)\tr$ and the vector of $\hat{\b w}_i$ is defined as follows depending on $z_i$, $x_i$ and the estimates from the first step:
\[
\hat{\b w}_i = \left\{
\begin{array}{ll}
(\hat{\pi}^*_{01}(0|\b v_i),\hat{\pi}^*_{01}(1|\b v_i),0,0)\tr, & \mbox{if } z_i=x_i=0,\\
(0,0,0,1)\tr, & \mbox{if } z_i=0,x_i=1,\\
(1,0,0,0)\tr, & \mbox{if } z_i=1,x_i=0,\\
(0,0,\hat{\pi}^*_{12}(1|\b v_i),\hat{\pi}^*_{12}(2|\b v_i))\tr, & \mbox{if } z_i=x_i=1.
\end{array}
\right. 
\]
Then we have the following first derivative:
\begin{eqnarray*}
\frac{\pa\ell_2(\b\be|\hat{\b\al}_0,\hat{\b\al}_2)}{\pa\b\be}&=&\sum_id_i\frac{\pa\ell_{2i}(\b\be|\hat{\b\al}_0,\hat{\b\al}_2)}{\pa\b\be},\\
\frac{\pa\ell_{2i}(\b\be|\hat{\b\al}_0,\hat{\b\al}_2)}{\pa\b\be}&=&\diag(\b a)\frac{\pa\ell_{2i}^*(\b\eta|\b\al_0,\b\al_2)}{\pa\b\eta},
\end{eqnarray*}
where $\b a=\diag(\b\eta)(\b 1-\b\eta)$, with $\b 1$ denoting a column vector of ones of suitable dimension. Similarly, with
\[
\frac{\pa^2\ell^*_2(\b\eta|\hat{\b\al}_0,\hat{\b\al}_2)}{\pa\b\eta\pa\b\eta\tr} = -\sum_id_i\left[\frac{y_{i2}}{(\hat{\b w}_i\tr\b\eta)^2}+\frac{1-y_{i2}}{(1-\hat{\b w}\tr_i\b\eta)^2}\right]\hat{\b w}_i\hat{\b w}_i\tr.
\]
we have that
\[
\frac{\pa^2\ell_2(\b\be|\hat{\b\al}_0,\hat{\b\al}_2)}{\pa\b\be\pa\b\be\tr} = \diag(\b a)\frac{\pa^2\ell_2^*(\b\eta|\hat{\b\al}_0,\hat{\b\al}_2)}{\pa\b\eta\pa\b\eta\tr}\diag(\b a)+\diag(\b b)\diag\left(\frac{\pa\ell_2^*(\b\eta|\hat{\b\al}_0,\hat{\b\al}_2)}{\pa\b\eta}\right),
\]
where $\b b=\diag(\b a)(\b 1-2\b\eta)$.

In order to compute standard errors for the parameter estimates, we use a sandwich formula for estimating the variance-covariance matrix of the overall estimator $\hat{\b\th}=(\b\al_0\tr,\b\al_1\tr,\b\be)\tr$. In particular, we have
\[
\hat{\b\Si}=\hat{\b H}^{-1}\hat{\b K}\hat{\b H}^{-1},
\]
where the matrices $\hat{\b H}$ and $\hat{\b K}$ are defined in Appendix.

Along the same lines as in \cite{bart:2010} we have performed a 
simulation study about the performance of the proposed estimator 
which we do not show here for reasons of space. Our
simulation study suggests good finite sample properties of the
estimator, also under more general assumptions than those formulated in Section 
\ref{sec:model}. Furthermore, it can be shown that if 
the control has the same effect on never-takers and
compliers, and the treatment has the same effect on compliers and
always-takers, the estimator $\hat{\b\be}$ is consistent as $n$ grows to infinity,
in symbols $\hat{\b\be}\convp\bar{\b\be}$, with $\bar{\b\be}=
(\bar{\be}_0,\bar{\be}_1,\bar{\be}_2,\bar{\be}_3)\tr$ denoting the true parameter
vector.

The result on existence and consistency of the
estimators is not ensured to hold when $\bar{\be}_0\neq\bar{\be}_1$ and/or
$\bar{\be}_2 \neq \bar{\be}_3$.
However, from the results of \cite{white:1982} on the maximum likelihood
estimation of misspecified models, it derives that $\hat{\b\be}\convp\b\be_*$,
where $\b\be_*$ is the supremum of $\E\{\ell_2(\b\be|\b\al_{0*},\b\al_{2*})/n\}$.
In the previous expression, $\b\al_{0*}$ and $\b\al_{2*}$
denote the limit in probability of $\hat{\b\al}_0$ and $\hat{\b\al}_2$, respectively.
We therefore expect $\b\be_*$
to be close to $\bar{\b\be}$ when $\bar{\be}_0$ is close to
$\bar{\be}_1$, $\bar{\be}_2$ is close 
to $\bar{\be}_3$, and $\pi(c|u,\b v)$ weakly depends on $u$. The same may
be said about the estimator $\hat{\de}$ of $\de$, whose limit in
probability is denoted by $\de_*$ and may be directly computed from
$\b\be_*$.

\section{Application to randomized study on coronary angiography after
myocardial infarction}\label{sec:application}
In this section we describe the application of the proposed estimator to the analysis of the data described in Section \ref{sec:data}. We recall that the proposed approach is
based on two steps: ({\em i}) estimation of the model for probability of being a never-taker, a complier, or an always-taker, and ({\em ii}) computation of the approximate conditional logistic estimator.

Regarding the first step, an important point is the selection of the
covariates to explain the non-compliance. 
In particular, we performed model choice by minimizing the Bayesian Information Criterion \citep[BIC, see][]{sch:78}, and finally selected two 
predictors (GI discretized using the quartiles and use of statines); see
also Section \ref{sec:data}. The
results from fitting this model are reported in Table
\ref{tab_ca_res_1} in terms of estimates of the parameters $\b\al_0$ and $\b\al_2$,
which are involved in expressions (\ref{eq:ass1}) and (\ref{eq:ass2}), and 
corresponding $t$-statistics and $p$-values.
For the categorical variable identifying the quartile of GI, we used the last quartile as reference category.

\begin{table}[!ht]\centering
{\small\begin{tabular}{ccccc}
\hline\hline        \vspace*{-0.4cm}                    \\
\multicolumn5c{\bf Parameter estimates for probability of being
  never-taker}\\
\hline
Estimator & Value & Std. Err. & $t$-statistic & $p$-value\\  \hline
$\hat{\al}_{00}$ (Intercept) &1.604  & 0.531  & 3.017 & 0.002 \\
$\hat{\al}_{01}$ (1st quartile GI) & -0.757 & 0.384 & -1.974 & 0.048 \\
$\hat{\al}_{02}$ (2nd quartile GI) & -0.886 & 0.368 & -2.406 & 0.016 \\
$\hat{\al}_{03}$ (3rd quartile GI) & -0.437 & 0.388 & -1.125 & 0.260 \\
$\hat{\al}_{04}$ (use of statin) & -0.985 & 0.438 & -2.247 & 0.025 \\\hline\\[1.5ex]
\hline\hline
\multicolumn5c{\bf Parameter estimates for probability of being always-taker}\\\hline
Estimator & Value & Std. Err. & $t$-statistic & $p$-value \\  \hline
$\hat{\al}_{20}$ (Intercept) & 1.454 & 0.597 & 2.436 & 0.015 \\
$\hat{\al}_{21}$ (1st quartile GI) & -0.565 & 0.444 & -1.274 & 0.202 \\
$\hat{\al}_{22}$ (2nd quartile GI) & -0.862 & 0.434 & -1.987 & 0.046 \\
$\hat{\al}_{23}$ (3rd quartile GI) & -0.459 & 0.449 & -1.023 & 0.306 \\
$\hat{\al}_{24}$ (use of statin) & -0.980 & 0.496 & -1.977 & 0.048 \\
\hline
\end{tabular}}
\caption{\em Estimates of compliance probability parameters for the proposed
  model, computed on the prompt coronary angiography data; predictors
  are quartiles of glicemic index (GI) and use of statines.}\label{tab_ca_res_1}
\end{table}

We observe a significant
non-compliance. The probabilities of being an always or a never taker 
are related in both cases with the GI and with use of statines. 
It can be seen that there is a significant lower probability of being
a always-taker in the second GI quartile with respect to the fourth,
while the other two quartiles are not statistically different from the
fourth. On the other hand, the probability of being a never taker steadily
increases with the GI category, with the third and fourth quartile not
being significantly different. The estimated effects of GI for always-takers
are explained considering that doctors may choose 
to assign to prompt CA even patients
randomized to the control (therefore making them always-takers) with
an abnormal GI (here, above the median or in the first quartile). 
Finally, the use of statines increases compliance in both
directions. This effect can be related to the fact that patients
using statines are better monitored and maybe already known to
doctors, and therefore an higher adherence to the experimental
settings is easier for these patients. 

Note that, even without
covariates, by the proposed method we can obtain an approximately
unbiased estimator of the causal effect (as seen by comparing $\hat\de$ with
$\hat\de^{(1)}$ in Table \ref{tab_ca_res_2}), but the use of covariates
allows to take into account part of the heterogeneity, therefore
decreasing the standard error of this estimate. 

\begin{table}[!ht]\centering
{\small\begin{tabular}{ccccc}
\hline\hline      
\multicolumn5c{\bf Estimates of the causal parameters}\\ 
\hline
Estimator & Value & Std. Err. & $t$-statistic & $p$-value \\  \hline
$\hat{\beta}_0$ &  2.158 & 0.361 &  5.973 & $<0.001$\\
$\hat{\beta}_1$ &  1.948 &0.677 & 2.878 &  0.004 \\
$\hat{\beta}_2$ & -0.072& 0.370 & -0.195 & 0.845 \\
$\hat{\beta}_3$ &  2.252 & 0.455 & 4.945 & $<0.001$ \\\hline\\[1.5ex]
\hline\hline
\multicolumn5c{\bf Estimates of the causal effect for compliers}\\
\hline
Estimator & Value & Std. Err. & $t$-statistic & $p$-value \\  \hline
$\hat{\de}$ (proposed method) &  -2.020 & 0.769 &  -2.625 & 0.009 \\ 
$\hat{\de}^{(1)}$ (proposed method) &  -1.938 & 0.929 &  -2.086 & 0.037 \\ 
$\hat{\de}^{(2)}$& -0.177 & 0.118 & -1.500 & 0.133 \\
$\hat{\de}^{(3)}$ & -0.513 & 0.119 & -4.311 & $<0.001$ \\
$\hat{\de}^{(4)}$ & -0.550 & 0.149 & -3.691 & $<0.001$\\
\hline
\end{tabular}}
\caption{\em Causal parameters for the proposed model estimated on the prompt coronary
  angiography data. Predictors are GI (discretized in quartiles) 
  and use of statines. In the bottom panel, $\hat\de$ is compared 
with the same estimate when covariates are not used ($\hat\de^{(1)}$)
and with competing estimators: $\hat{\de}^{(2)}$ standard conditional estimator based 
on received treatment ($X$);  $\hat{\de}^{(3)}$ standard
conditional estimator based on assigned treatment ($Z$, Intention to
Treat analysis); $\hat{\de}^{(4)}$ standard conditional estimator based on the
assigned and complied treatment (Per Protocol analysis)}\label{tab_ca_res_2}
\end{table} 

In Table \ref{tab_ca_res_2} we report estimates of causal parameters, 
and compare them with four other estimators. The first 
(denoted by $\hat{\de}^{(1)}$) is based on 
our proposed approach in which no covariates are used to predict
compliance. The other three estimators (denoted by $\hat{\de}^{(2)}$, $\hat{\de}^{(3)}$,
and $\hat{\de}^{(4)}$, respectively)
are based on conditional logistic regression on the received
treatment, an Intention to Treat and a Per Protocol analysis. The last
two are based on the assigned treatment regardless of the actually
received treatment, and on patients actually receiving the assigned
treatment, respectively. From the upper panel we 
can see that the control has approximately the same effect on never
takers and on compliers (with a log-odds of about 2). 
The treatment seems to have no effect on compliers, while on 
always-takers we once again obtain a log-odds of about 2. 
We therefore can say that 
({\em i}) lack of a prompt CA, regardless of whether it was assigned
or as a result of non-compliance, may increase the risk of recurrence
and ({\em ii}) if a patient who was assigned to the control group undergoes
prompt CA, this is likely due to a possibly bad (even life
threatening) condition, hence the high risk of recurrent events even
under the treatment. A consequence is that bias 
with ITT and PP estimators arise mostly due to always
takers. In fact, the effect of the control is approximately the same
on never-takers and compliers ($\beta_0\approx\beta_1$); on 
the other hand, there is a strong
difference of the effect of treatment as estimated on compliers and
always-takers ($\beta_2 \neq \beta_3$). 

Always-takers in this example can be 
expected to experience the event even after the treatment. 
Ignoring this fact
will make the two groups artificially more similar, as testified by
the estimates $\hat\de$$^{(2)}$, $\hat\de$$^{(3)}$, and $\hat\de$$^{(4)}$. 
In fact, our most important estimate is $\hat\de$, which is approximately
$-2$. When our final estimate is compared with 
$\hat\de$$^{(2)}$, $\hat\de$$^{(3)}$, and $\hat\de$$^{(4)}$ 
we find that those are 
at most only half our causal estimate. The estimate of the causal
parameter based on the received treatment ($\hat\de$$^{(2)}$) is not
even significant. 
Standard fits in this example may lead to grossly
underestimate the effect of a prompt CA. 

%
\section{Discussion}\label{sec:conclusions}
An approach has been introduced to estimate the causal effect of a
treatment over control on the basis of a two-arm experimental study
with possible non-compliance. The approach is applicable when the
effect of the treatment is measured by a binary response variable
observed before and after the treatment. It relies on a causal model
formulated on the basis of latent variables for the effect of
unobservable covariates at both occasions and to account for the
difference between compliers and non-compliers in terms of reaction
to control and treatment. The parameters of the model are estimated
by a pseudo conditional likelihood approach based on an approximated
version of the conditional probability of the two response variables
given their sum. The causal model and the proposed estimator extend the
model and the estimator of \cite{bart:2010} to the case in which non-compliance
may also happen in the control arm.

The method is applied to the analysis of data coming from a study
on the effect of prompt coronary angiography in myocardial
infarction. The application shows that prompt coronary angiography in patients with
myocardial infarction may significantly decrease the risk of other
events within the next two years, with a log-odds of about -2. On the other hand,
estimates of this log-odds ratio obtained by the standard logistic
approach are considerably closer to 0.

One of the basic assumptions on which the approach relies is that a
subject is assigned to the control arm or to the treatment arm with
a probability depending only on the observable covariates and not on
the pre-treatment response variable. Indeed, we could relax this
assumption, but we would have much more complex expressions for the
conditional probability of the response variables given their sum.

As a final comment we remark that  
we only considered the case of repeated response
variables having a binary nature. However, the approach may be
directly extended to the case of response variables having a different
nature (e.g. counting), provided that the conditional distribution
of these variables belongs to the natural exponential family and the
causal effect is measured on a scale defined according to the
canonical link function for the adopted distribution \citep{mcc:nel:89}.

\subsection*{Appendix: Matrices involved in the sandwich estimator for the variance of the estimator}
We have that
\[
\hat{\b H}=\pmatrix{
{\displaystyle\frac{\pa\ell_1(\hat{\b\al}_0,\hat{\b\al}_2)}{\pa\b\al_0\pa\b\al_0\tr}} & 
{\displaystyle\frac{\pa\ell_1(\b\al_0,\b\al_2)}{\pa\b\al_0\pa\b\al_2\tr}} & \b O\cr
{\displaystyle\frac{\pa\ell_1(\b\al_0,\b\al_2)}{\pa\b\al_2\pa\b\al_0\tr}} & 
{\displaystyle\frac{\pa\ell_1(\b\al_0,\b\al_2)}{\pa\b\al_2\pa\b\al_2\tr}} & \b O\cr
{\displaystyle\frac{\pa\ell_2(\b\be|\hat{\b\al}_0,\hat{\b\al}_2)}{\pa\b\be\pa\b\al_0\tr}} & 
{\displaystyle\frac{\pa\ell_2(\b\be|\hat{\b\al}_0,\hat{\b\al}_2)}{\pa\b\be\pa\b\al_2\tr}} & 
{\displaystyle\frac{\pa\ell_2(\b\be|\hat{\b\al}_0,\hat{\b\al}_2)}{\pa\b\be\pa\b\be\tr}} 
}
\]
and
\[
\hat{\b K} = \sum_i \pmatrix{{\displaystyle\frac{\pa\ell_{1i}(\hat{\b\al_0})}{\pa\b\al_0}} \cr 
{\displaystyle\frac{\pa\ell_{1i}(\hat{\b\al_2})}{\pa\b\al_2}}\cr
{\displaystyle\frac{\pa\ell_{2i}(\hat{\b\be}|\hat{\b\al}_0,\hat{\b\al}_2)}{\pa\b\be}}}
\pmatrix{{\displaystyle\frac{\pa\ell_{1i}(\hat{\b\al_0})}{\pa\b\al_0\tr}} & {\displaystyle\frac{\pa\ell_{1i}(\hat{\b\al_2})}{\pa\b\al_2\tr}} &
{\displaystyle\frac{\pa\ell_{2i}(\hat{\b\be}|\hat{\b\al}_0,\hat{\b\al}_2)}{\pa\b\be\tr}}}.
\]
In the above expressions, $\b O$ denotes a matrix of zeros of suitable dimension. Moreover, all the derivatives have been defined, with the exception of the derivative of $\ell_2(\hat{\b\be}|\hat{\b\al}_0,\hat{\b\al}_2)$ with respect to $\b\al_0$ (or $\b\al_2$) and $\b\be$. In particular, we have that:
\[
\frac{\pa^2\ell_2^*(\b\be|\hat{\b\al}_0,\hat{\b\al}_2)}{\pa\b\be\pa\b\al_c\tr} = \diag(\hat{\b a})\sum_id_i\left(\frac{y_{i2}}{\hat{\b w}_i\tr\b\eta}-\frac{1-y_{i2}}{1-\hat{\b w}\tr_i\b\eta}\right)\frac{\pa\hat{\b w}_i}{\pa\b\al_c\tr},\quad c = 0,2,
\]
where
\[
\frac{\pa\hat{\b w}_i}{\pa\b\al_0\tr} = \left\{
\begin{array}{ll}
(\hat{\pi}^*_{01}(0|\b v_i)\hat{\pi}^*_{01}(1|\b v_i),-\hat{\pi}^*_{01}(0|\b v_i)\hat{\pi}^*_{01}(1|\b v_i),0,0)\tr\b g(\b v_i), & z_i=x_i=0,\\
\b O, & \mbox{otherwise},
\end{array}
\right. 
\]
and
\[
\frac{\pa\hat{\b w}_i}{\pa\b\al_0\tr} = \left\{
\begin{array}{ll}
(0,0,-\hat{\pi}^*_{12}(1|\b v_i)\hat{\pi}^*_{01}(1|\b v_i),\hat{\pi}^*_{12}(1|\b v_i)\hat{\pi}^*_{12}(2|\b v_i))\tr\b g(\b v_i), & z_i=x_i=1,\\
\b O, & \mbox{otherwise}.
\end{array}
\right. 
\]

\bibliography{biblio}
\bibliographystyle{apalike}
\end{document}